\documentclass{amsart}

\usepackage[foot]{amsaddr}

\usepackage{tikz}
\usepackage{graphicx}
\usepackage{enumerate,amsfonts}

\newcommand{\Z}{\mathbb{Z}}
\newcommand{\Q}{\mathbb{Q}}
\newcommand{\C}{\mathbb{C}}
\renewcommand{\a}{\alpha}
\renewcommand{\b}{\beta}
\renewcommand{\c}{\gamma}
\renewcommand{\t}{\theta}
\newcommand{\Norm}{\mathcal{N}}
\newcommand{\Ia}{\mathfrak{a}}
\newcommand{\Ip}{\mathfrak{p}}

\newtheorem{thm}{Theorem}[section]

\newtheorem{cor}[thm]{Corollary}
\newtheorem{prop}[thm]{Proposition}
\newtheorem{rmk}[thm]{Remark}
\newtheorem{exmp}{Example}
\theoremstyle{definition}
\newtheorem{defn}[thm]{Definition}

\title{First-degree prime ideals of biquadratic fields dividing prescribed principal ideals}

\author{Giordano Santilli$^1$}
\email{giordano.santilli@unitn.it}
\author{Daniele Taufer$^1$}
\email{daniele.taufer@gmail.it}

\address{%
 $^1$University of Trento}

\begin{document}
\maketitle

\begin{abstract}We describe first-degree prime ideals of biquadratic extensions in terms of first-degree prime ideals of two underlying quadratic fields. The identification of the prime divisors is given by numerical conditions involving their ideal norms. The correspondence between these ideals in the larger ring and those in the smaller ones extends to the divisibility of special-shaped principal ideals in their respective rings, with some exceptions that we explicitly characterize.
\end{abstract}
\section{Introduction}
Biquadratic fields are numerical fields that have been studied extensively \cite{BaeYue,Dwilewicz,Mamoona,Ouyang,Williams, Yue} and are currently in the spotlight as they provide examples of non-principal Euclidean ideal classes \cite{ChaMut,HuLi}.
They are defined as numerical fields whose Galois group is the Klein group and they can be obtained by the compositum of two quadratic fields, which is the construction we adopt in the present work.

In a general number field $\Q(\c)$, the special subring  $\Z[\c]$ of the ring of integers plays an important role in applications, since their elements have a natural representation as integer-valued polynomials. Among its main properties, identifying its first-degree prime ideals is relatively easy as well as deciding whether they divide a given principal ideal \cite{BLP}.
This feature is crucial for the effectiveness of the General Number Field Sieve, which is the most efficient known algorithm to factorize large integers \cite{LLM,LLMP,Kleinjung}.

In this paper, we start from a biquadratic field $\Q(\c)$ and we investigate the relation between the first-degree prime ideals in $\Z[\c]$ and the first-degree prime ideals of $\Z[\a]$ and $\Z[\b]$, with $\Q(\a)$ and $\Q(\b)$ two underlying quadratic fields.

\begin{center}
\begin{tikzpicture}[node distance=2cm]
 \node at (0,0) (Q) {$\Q$};
 \node at (-1,1) (Qt1){$\Q(\a)$}; 
 \node at (1,1) (Qt2) {$\Q(\b)$};
 \node at (0,2) (Qt) {$\Q(\c) = \Q(\a +\b)$};
 \draw (Q) -- (Qt1) node [midway, below, xshift = -0.15cm, yshift = 0.1cm] (TextNode) {\footnotesize 2};
 \draw (Qt1) -- (Qt) node [midway, above, xshift = -0.15cm, yshift = -0.1cm] (TextNode) {\footnotesize 2};
 \draw (Q) -- (Qt2) node [midway, below, xshift = 0.15cm, yshift = 0.1cm] (TextNode) {\footnotesize 2};
 \draw (Qt2) -- (Qt) node [midway, above, xshift = 0.15cm, yshift = -0.1cm] (TextNode) {\footnotesize 2};
\end{tikzpicture}
\end{center}

In Section 2, we recall the results we need on the structure of first-degree prime ideals.
In Section 3, we provide an explicit relation between such ideals in $\Z[\c]$ and those in $\Z[\a]$ and $\Z[\b]$, which depends on their ideal norms.
In Section 4, we extend this relation to the divisibility of principal ideals in their respective rings, with some exceptions that we explicitly highlight.
In Section 5, we hint at applications and further research.

\section{Notation and preliminaries}
Let $\Q(\a)$ and $\Q(\b)$ be two distinct quadratic fields, i.e. $a = \a^2$ and $b = \b^2$ have distinct non-trivial square-free parts. We may assume these number fields to be generated by the polynomials

\begin{equation*}
    f_a(x) = x^2 - a, \quad \quad f_b(x) = x^2 - b,
\end{equation*}
where $a,b \in \Z$.
It is well-known \cite{Williams} that the biquadratic extension they generate is $\Q(\c)$ with $\c = \a + \b$, whose minimal polynomial is

\begin{equation*}
    f_c(x) = x^4 - 2(a+b)x^2 + (a-b)^2.
\end{equation*}

In this setting we focus on ideals in the order $\Z[\t]$, which has been deeply studied in \cite{BLP}.

\begin{defn} 
Let $\t$ be an algebraic integer and $\Ia$ be a non-zero ideal of $\Z[\t]$. The norm of $\Ia$ is

\begin{equation*}
    \Norm(\Ia) = [\Z[\t] : \Ia].
\end{equation*}
Moreover, a non-zero prime ideal $\Ip$ of $\Z[\t]$ is called a first-degree prime ideal if $\Norm(\Ip)$ is a prime integer.
\end{defn}

The following theorem characterizes the first-degree prime ideals of $\Z[\t]$.

\begin{thm}\textnormal{(\cite[p.57]{BLP})}
Let $f \in \Z[\t]$ be an irreducible monic polynomial and $\t \in \C$ one of its roots. Then, for every positive prime $p$ there is a bijection between

\begin{equation*}
    \{(r,p) \ | \ r \in \Z / p\Z \text{ such that } f(r) \equiv 0 \bmod p\}
\end{equation*}
and

\begin{equation*}
    \{ \Ip \ | \ \Ip \in \textnormal{Spec}\ \Z[\t] \text{ such that } \Norm(\Ip) = p \}.
\end{equation*}
\end{thm}

In the above bijections an ideal $\Ip$ corresponds to $(r,p)$ if it is the kernel of the evaluation-in-$r$ ring morphism

\begin{equation*}
    \pi_{\t}: \Z[\t] \to \Z/p\Z, \quad a+b\t \mapsto a+br \bmod p.
\end{equation*}
In the following sections we will use these bijections as identifications.

It is a standard fact \cite{Dedekind} that non-zero ideals factor uniquely into primes in the whole ring of integers, however such a decomposition has been generalized even inside non-maximal orders $\Z[\t]$ in \cite[Proposition 7.1]{BLP}.
The relation between these two types of factorization has been precisely characterized in \cite[Proposition 7.2]{BLP}, and it founds a relevant employment for factoring in $\Z[\t]$ special principal ideals \cite[Corollary 5.5]{BLP}, whose factorization only consists of first-degree prime ideals in the same order.

\section{First-degree prime ideals of biquadratic extensions}

The following theorem exhibits how to construct first-degree prime ideals of a given norm in a biquadratic field by knowing first-degree prime ideals of the same norm in two of its quadratic subfields.
\begin{thm} \label{thm1}
Let $(r,p)$ be a first-degree prime ideal of $\Z[\a]$ and $(s,p)$ a first-degree prime ideal of $\Z[\b]$. Then $(r+s,p)$ is a first-degree prime ideal of $\Z[\c]$.
\end{thm}
\begin{proof}
By hypothesis we have 

\begin{align*}
  f_a(r) = r^2-a \equiv 0 \bmod{p}, \\
  f_b(s) = s^2-b \equiv 0 \bmod{p}. 
\end{align*}
By plugging $r+s$ into $f_c$ we get

 \begin{align*}
  f_c(r+s)&=(r+s)^4-2(a+b)(r+s)^2+(a-b)^2 \\
  &\equiv (r+s)^4 - 2(r^2+s^2)(r+s)^2 + (r^2-s^2)^2 \bmod{p} \\
  &\equiv (r+s)^2 \big( (r+s)^2 - 2(r^2+s^2) + (r-s)^2 \big) \bmod{p} \\
  &\equiv 0 \bmod{p}.
 \end{align*}
 Therefore, $(r+s,p)$ is a first-degree prime ideal of $\Z[\c]$.
\end{proof}

We will refer to $(r+s, p) \subseteq \Z[\c]$ as the \emph{combination} of the ideals $(r, p) \subseteq \Z[\a]$ and $(s, p) \subseteq \Z[\b]$.
Now we prove that the ideals constructed as combinations are almost all the first-degree prime ideals of $\Z[\c]$.

\begin{thm} \label{thm2}
Let $(t,p)$ be a first-degree prime ideal of $\Z[\c]$. If either $p = 2$ or $t \not \equiv 0 \bmod p$ then there exists a unique pair $r,s \in \Z / p\Z$ such that $t \equiv r + s \bmod p$ and $(r,p), (s,p)$ are first-degree prime ideals of $\Z[\a]$ and $\Z[\b]$, respectively.
\end{thm}
\begin{proof}
We treat separately the cases $p = 2$ and $p$ odd, explicitly exhibiting such a pair $(r,s)$ in both cases.
\begin{itemize}
 \item \textbf{Case: $p = 2$}.
 
 Since in $\Z / 2\Z$ every element is equal to its square, the only choice of $r,s \in \Z / 2\Z$ satisfying
 
 \begin{equation*}
     \begin{cases}
     r^2 - a \equiv 0 \bmod 2,\\
     s^2 - b \equiv 0 \bmod 2,
     \end{cases}
 \end{equation*}
 is $(r,s) = (a,b)$, which also satisfies $t = r+s$ because
 
\begin{equation*}
    0 \equiv f_c(t) \equiv t^4-2(a+b)t^2+(a-b)^2 \equiv t + a + b \bmod 2.
\end{equation*}
 
 \item \textbf{Case: $p \neq 2$ and $t \neq 0$}.
 
 In this case $2t$ is invertible in $\Z/p\Z$, then we can define $r_1 = \frac{t^2+a-b}{2t}$.
 We notice that from $f_c(t) \equiv 0 \bmod p$ we have
 
 \begin{align*}
    r_1^2 = \left( \frac{t^2+a-b}{2t} \right)^2 &= \frac{t^4+2at^2-2bt^2+a^2-2ab+b^2}{4t^2}\\
    &= \frac{f_c(t) + 4at^2}{4t^2} \equiv a \bmod p.
 \end{align*}
 Thus, $r_1$ is a square root of $a$ modulo $p$ and since $\Z / p\Z$ is a finite field there are at most two solutions to $r^2-a \equiv 0 \bmod p$, hence these are
 
 \begin{equation*}
     r_1 = \frac{t^2+a-b}{2t}, \quad r_2 = -r_1 = -\frac{t^2+a-b}{2t}.
 \end{equation*} 
 Similarly, there are only two possible values for $s$, which are
 
 \begin{equation*}
     s_1 = \frac{t^2-a+b}{2t}, \quad s_2 = -s_1 = -\frac{t^2-a+b}{2t}.
 \end{equation*} 
 It is easy to verify that $r_1 + s_1 = t$ and we now prove that $(r_1,s_1)$ is in fact the unique choice for such a pair $(r,s)$ in order to satisfy $r + s = t$.
 
 First, we notice that $(r,s) = (r_2, s_2)$ is not a possible option, since in this case $r_2 + s_2 = -t$ but $-t \not \equiv t \bmod p$ since $p \neq 2$ and $t \not \equiv 0 \bmod p$.
 
 To conclude the proof we show that $(r_1, s_2)$ may be a suitable choice only when $s_1 = s_2 = 0$, therefore $(r_1, s_2) = (r_1, s_1)$.
 Suppose that $r_1 + s_2 = t$,  which means $a - b = t^2$. Then
 
 \begin{equation*}
     0 \equiv f_c(t) \equiv 2t^2(t^2-a-b) \bmod p.
 \end{equation*}
 Since $2t^2 \not \equiv 0 \bmod p$ we get $a + b \equiv t^2 \bmod p$, then $a \equiv t^2 \bmod p$ and $b \equiv 0 \bmod p$. This proves that $s_1 = s_2 = 0$.
 
 The same argument shows that $(r_2,s_1)$ is a valid pair only if $r_1 = r_2 = 0$ so $(r_2,s_1) = (r_1,s_1)$.
 In conclusion, there is only one working pair, that is $(r,s) = (r_1, s_1)$. \qedhere
\end{itemize} 
\end{proof}

The uniqueness part of Theorem \ref{thm2} states that any ideal $(t,p)$ of $\Z[\c]$ with $t \neq 0$ may be determined without repetitions from first-degree prime ideals of two underlying quadratic fields.
The only ideals left are those of the form $(0,p)$ for $p \neq 2$, which are examined in the following proposition.

\begin{prop} \label{prop1}
Let $(0,p)$ be a first-degree prime ideal of $\Z[\c]$ and $r \in \Z / p\Z$, then the following are equivalent:
\begin{itemize}
    \item $f_a(r) \equiv 0 \bmod p$,
    \item $f_b(r) \equiv 0 \bmod p$,
    \item $(r,p)$ and $(-r,p)$ are first-degree prime ideals of $\Z[\a]$,
    \item $(r,p)$ and $(-r,p)$ are first-degree prime ideals of $\Z[\b]$.
\end{itemize}
\end{prop}
\begin{proof}
From $f_c(0) \equiv 0 \bmod p$ we get

\begin{equation*}
    (a - b)^2 \equiv 0 \bmod p \implies a \equiv b \bmod p.
\end{equation*}
Hence, $f_a \equiv f_b \bmod p$, therefore $r$ is a root of $f_a$ modulo $p$ if and only if the same holds for $f_b$. Moreover, if $f_a(r) \equiv 0 \bmod p$ also $f_a(-r) \equiv 0 \bmod p$, implying that $(\pm r, p)$ are first-degree prime ideals of $\Z[\a]$, while the converse is trivial.
\end{proof}

\begin{rmk}
The above proposition is trivial for $p = 2$. In fact, by Theorem \ref{thm2} if $(0,2)$ is a first-degree prime ideal of $\Z[\c]$ then all the above equivalent conditions are satisfied for $r = a$.
\end{rmk}

According to Proposition \ref{prop1}, one of the following situations takes place, depending on the number $\nu$ of roots of $f_a$ modulo $p$:
\begin{itemize}
    \item[\footnotesize $\nu = 0$:] $(0,p) \subseteq \Z[\c]$ cannot be found as a combination of first-degree prime ideals of $\Z[\a]$ and $\Z[\b]$.
    \item[\footnotesize $\nu = 1$:] $(0,p) \subseteq \Z[\c]$ is the combination of $(0,p) \subseteq \Z[\a]$ and $(0,p) \subseteq \Z[\b]$.
    \item[\footnotesize $\nu = 2$:] $(0,p) \subseteq \Z[\c]$ is determined by two different combinations of first-degree prime ideals of $\Z[\a]$ and $\Z[\b]$.
\end{itemize}

In the following example, we see that all the above cases may actually occur.

\begin{exmp}
Let $f_a = x^2 - 50$ and $f_b = x^2 - 155$ generate the quadratic fields $\Q(\a)$ and $\Q(\b)$, so that the composite biquadratic field $\Q(\c)$ is generated by the polynomial $f_c = x^4 - 410x^2 + 11025$.

The unique first-degree prime ideal in $\Z[\c]$ with norm $p = 3$ is $(0,3)$, but there are no such ideals neither in $\Z[\a]$ nor in $\Z[\b]$, then $(0,3)$ cannot be a combination of any of them.

The unique first-degree prime ideal of norm $p = 5$ in $\Z[\c]$ is $(0,5)$, which is determined uniquely as a combination of the ideals $(0,5)$ in $\Z[\a]$ and $(0,5)$ in $\Z[\b]$.

There are three first-degree prime ideals of norm $p = 7$ in $\Z[\c]$: $(0,7), (2,7)$ and $(5,7)$.
The first-degree prime ideals of the same norm for both $\Z[\a]$ and $\Z[\b]$ are $(1,7)$ and $(6,7)$.
As prescribed by Theorem \ref{thm2} we observe that $(2,7)$ and $(5,7)$ are uniquely determined by the combinations of $\big((1,7), (1,7)\big)$ and $\big((6,7), (6,7)\big)$, whereas $(0,7)$ arises from both the combinations of $\big((1,7), (6,7)\big)$ and $\big((6,7), (1,7)\big)$.
\end{exmp}

\section{Division of prescribed principal ideals}

In this section we consider a special family of principal ideals of $\Z[\c]$ and we study its first-degree prime divisors in terms of first-degree prime ideals dividing its intersections with $\Z[\a]$ and $\Z[\b]$. This particular class of ideals is of great interest, since they are the ones usually employed in applications \cite{BLP}. 
We begin by characterizing these intersections, as in the following proposition.

\begin{prop} \label{PropPrincipal}
\label{IntersectionProp} 
Let $n$ and $m \neq 0$ be coprime integers and let $I = \langle n + m\c \rangle \subseteq \Z[\c]$. Then $I \cap \Z[\a]$ is a principal ideal of $\Z[\a]$ generated by 

\begin{equation*}
I \cap \Z[\a] = \left\langle (n + m\a + m\b)(n + m\a - m\b) \right\rangle.
\end{equation*}
\end{prop}
\begin{proof}
$\supseteq$ ) The generator $(n + m\a + m\b)(n + m\a - m\b)$ is an element of $I$ and is equal to $(n+m\a)^2 - (m\b)^2 = n^2+m^2(a-b) + 2nm\a$, which belongs to $\Z[\a]$.\\
$\subseteq$ ) Any element $x \in I$ may be written as 

\begin{equation*}
    x = (n+m\a+m\b)(\lambda_0 + \lambda_1 \a + \lambda_2 \b + \lambda_3 \a \b),
\end{equation*}
for some $\lambda_0, \dots, \lambda_3 \in \Z$. Since $\{1, \a, \b, \a\b\}$ is a $\Q$-basis of $\Z[\c]$ \cite{Williams}, then $x \in \Z[\a]$ if and only if its coefficients of $\b$ and $\a\b$ are zero, which amounts to

\begin{equation}
\label{eq:intersection1}
    m\lambda_0 + n\lambda_2 + ma\lambda_3 = 0, \quad  m\lambda_1 + m\lambda_2 + n\lambda_3 = 0.
\end{equation}
From Equation \eqref{eq:intersection1} we get

\begin{align*}
    \lambda_0 + \lambda_1 \a + \lambda_2 \b + \lambda_3 \a \b 
    &= -\left(\frac{n}{m}\lambda_2 + a\lambda_3\right) - \a \left(\lambda_2 + \frac{n}{m}\lambda_3\right) + m\b \left(\frac{\lambda_2 + \lambda_3 \a}{m}\right)\\
    &= (n + m\a - m\b)\left(-\frac{\lambda_2 + \lambda_3\a}{m}\right).
\end{align*}
Moreover since $(m,n)=1$, then Equation \eqref{eq:intersection1} also implies that $m$ divides both $\lambda_2$ and $\lambda_3$.  
Hence, we conclude that any $x \in I \cap \Z[\a]$ belongs to the ideal generated in $\Z[\a]$ by $(n + m\a + m\b)(n + m\a - m\b)$.
\end{proof}

With the following theorems we prove that divisibility is stable under combination except for an exceptional case.

\begin{thm}
 Let $n$ and $m \neq 0$ be coprime integers and $I = \langle n + m\c \rangle$ be a principal ideal of $\Z[\c]$. Let us assume that there are $(r, p)$ first-degree prime ideal of $\Z[\a]$ dividing $I_a = I \cap \Z[\a]$ and $(s, p)$ first-degree prime ideal of $\Z[\b]$ dividing $I_b = I \cap \Z[\b]$.
 Then $(r + s, p)$ is a first-degree prime ideal of $\Z[\c]$ dividing $I$ unless the following conditions simultaneously hold:
 
 \begin{equation*}
     p \neq 2, \quad n \equiv 0 \bmod p, \quad r + s \not \equiv 0 \bmod p.
 \end{equation*}
\end{thm}

\begin{proof}
 By Theorem \ref{thm1} $(r + s, p)$ is a first-degree prime ideal of $\Z[\c]$, then it is sufficient to show that under the aforementioned conditions we obtain
 
 \begin{equation} \label{eqtp}
 n+m(r+s) \equiv 0 \bmod p,
 \end{equation}
which proves that $I = \langle n + m\c \rangle \subseteq \ker \pi_{\c} = (r + s, p)$, hence $(r + s, p) | I$.

\begin{itemize}
 \item \textbf{Case:} $p = 2$.
 
 By Proposition \ref{PropPrincipal} a generator of $I_a$ in $\Z[\a]$ is $g = n^2 + m^2(a-b) + 2nm\a$, therefore
 
 \begin{equation*}
    \pi_{\a}(g) = n^2 + m^2(a-b) + 2nmr \equiv n + m(a^2+b^2) \equiv n + m(r+s) \bmod 2.
 \end{equation*}
 Thus, either \eqref{eqtp} is satisfied or $\pi_{\a}(g) \equiv 1 \bmod 2$, which implies that there are no first-degree prime ideals $(r,2)$ dividing $I_a$.

 \item \textbf{Case}: $p \neq 2$ and $n \not \equiv 0 \bmod p$.
 
By Proposition \ref{PropPrincipal} we have 

\begin{equation*}
\begin{cases}
I_a = \langle n^2+m^2(a-b)+2nm\a \rangle \subseteq \Z[\a], \\
I_b = \langle n^2+m^2(b-a)+2nm\b \rangle \subseteq \Z[\b].
\end{cases}
\end{equation*}
Since $\ker \pi_{\a} = (r,p) \mid I_a$ and $\ker \pi_{\b} = (s,p) \mid I_b$, then

\begin{equation*}
 \begin{cases}
  n^2+m^2(a-b)+2nmr \equiv 0 \bmod p, \\
  n^2+m^2(b-a)+2nms \equiv 0 \bmod p.
 \end{cases}
\end{equation*}
Summing the above relations we get

\begin{equation*}
 2n[n+m(r+s)] \equiv 0 \bmod p.
\end{equation*}
Since $2n \not \equiv 0 \bmod p$ this implies \eqref{eqtp}.

\item \textbf{Case}: $p \neq 2$, $n \equiv 0 \bmod p$ and $r + s \equiv 0 \bmod p$.
 
 In this case \eqref{eqtp} is trivially satisfied.
\end{itemize}
Thus, we conclude that $(r + s, p)$ is a first-degree prime ideal of $\Z[\c]$ dividing $I$ except for the case $p \neq 2$, $n \equiv 0 \bmod p$ and $r + s \not \equiv 0 \bmod p$.
\end{proof}

The next example shows that in the exceptional case mentioned above ideal combination may not maintain divisibility.

\begin{exmp}
Let $f_a = x^2 + 4$ and $f_b = x^2 - 6$ generate the quadratic fields $\Q(\a)$ and $\Q(\b)$, so that the composite biquadratic field $\Q(\c)$ is generated by the polynomial $f_c = x^4 - 4x^2 + 100$.

The first-degree prime ideals of $\Z[\c]$ with norm $p = 5$ are $(0,5), (2,5)$ and $(3,5)$, while those of $\Z[\a]$ and $\Z[\b]$ are $(1,5)$ and $(4,5)$.

Let $I$ be the principal ideal $\langle 5 + \c \rangle \subseteq \Z[\c]$. By Proposition \ref{PropPrincipal} we have 

\begin{equation*}
I_a = \langle 15 + 10\a \rangle \subseteq \Z[\a], \quad I_b = \langle 35 + 10\b \rangle \subseteq \Z[\b].
\end{equation*}
It is easy to see that both $(1,5)$ and $(4,5)$ divide $I_a$ and $I_b$.
Besides, the combination of $(1,5)$ with $(4,5)$ is $(0,5)$, which divides $I$.
However, the other options are exceptional, since the combination between $(1,5)$ and $(1,5)$ is $(2,5)$, which does not divide $I$.
The same holds for $(3,5)$, which is the combination of $(4,5)$ with $(4,5)$.
\end{exmp}

On the other hand, whenever a first-degree prime ideal of $\Z[\c]$ dividing a given principal ideal $I$ is obtained as a combination of two first-degree prime ideals, they divide the intersections of $I$ with $\Z[\a]$ and $\Z[\b]$.

\begin{thm} \label{thml}
 Let $n$ and $m\neq 0$ be coprime integers, $I = \langle n + m\c \rangle \subseteq \Z[\c]$ and let $(t,p)$ be a first-degree prime ideal dividing $I$. If there exist first-degree prime ideals $(r,p) \subseteq \Z[\a]$ and $(s,p) \subseteq \Z[\b]$ such that $r + s \equiv t \bmod{p}$, then $(r,p)$ divides $I_a = I \cap \Z[\a]$ and $(s,p)$ divides $I_b = I \cap \Z[\b]$.
\end{thm}
\begin{proof}
If these ideals exist, from $(t, p) \mid I$ we get

\begin{equation*}
    0 \equiv n + mt \equiv n + mr + ms \bmod p.
\end{equation*}
Let $g_a = n^2 + m^2(a-b) + 2nm\a$ a generator of $I_a$ and $g_b = n^2 + m^2(b-a) + 2nm\b$ a generator of $I_b$. From the above equation we have

\begin{align*}
    \pi_{\a}(g_a) &= n^2 + m^2(a-b) + 2nmr \equiv n^2 + m^2(a-b) + 2n(-n-ms) \bmod p \\
    &\equiv -n^2 - m^2(-a+b) - 2nms \equiv - \pi_{\b}(g_b) \bmod p.
\end{align*}
Hence, since $\pi_{\a}(g_a)$ vanishes if and only if $\pi_{\b}(g_b)$ does, it is sufficient to show that $(r,p) \mid I_a$. Substituting $n = - mr - ms$, $r^2 \equiv a \bmod p$ and $s^2 \equiv b \bmod p$ we get

\begin{equation*}
    \pi_{\a}(g_a) \equiv (-mr - ms)^2 + m^2(r^2-s^2) + 2(- mr - ms)mr \equiv 0 \bmod p,
\end{equation*}
therefore $(r,p) \mid I_a$.
\end{proof}

The following corollary enhances the previous result in the generic case.

\begin{cor}
Let $n$ and $m\neq 0$ be coprime integers, $I = \langle n + m\c \rangle \subseteq \Z[\c]$ and let $(t,p)$ be a first-degree prime ideal dividing $I$, with $t \neq 0$ if $p \neq 2$. Then there exist two unique first-degree prime ideals $(r,p) \subseteq \Z[\a]$ and $(s,p) \subseteq \Z[\b]$ such that $(r,p)$ divides $I \cap \Z[\a]$, $(s,p)$ divides $I \cap \Z[\b]$ and $r + s \equiv t \bmod{p}$.
\end{cor}
\begin{proof}
It follows immediately from Theorem \ref{thml} and Theorem \ref{thm2}.
\end{proof}

\begin{rmk}
For completeness, we comment on a case we have not discussed, which is $m = 0$ so that $I = \langle n \rangle \subseteq \Z[\c]$.
Both the intersections of this principal ideal with $\Z[\a]$ and $\Z[\b]$ are equal to $\langle n \rangle$ and for every $t \in \Z / p\Z$ we have 

\begin{equation*}
    (t,p) | \langle n \rangle \quad \iff \quad n \equiv 0 \bmod p.
\end{equation*}
Thus, all or none the pairs $\{(t,p) \ | \ f_c(t) \equiv 0 \bmod p\}$ divide $\langle n \rangle \subset \Z[\c]$, so it is still true that the combination preserves divisibility, but in general not uniquely.
\end{rmk}

\section{Applications and further work}

In this work we show how two first-degree prime ideals in quadratic extensions may be combined to obtain a first-degree prime ideal in the corresponding biquadratic extension lying over them. In addition, this correspondence is proved to preserve the division of prescribed first-degree prime ideals, except for some sporadic, though well-determined, cases.

Nonetheless, further computations suggest that our results might also be extended to more general number fields, possibly requiring additional hypotheses.
Such a generalization could be repeatedly applied in order to characterize first-degree prime ideals of a given norm in large extensions: the (at most) $d_1d_2$ first-degree prime ideals of a composite extension may be seen as combinations of $d_1 + d_2$ such ideals into smaller subrings, which are much more convenient to be stored and managed.

Among other applications, first-degree prime ideals are employed in the General Number Field Sieve, where a large number of them is needed to factorize some principal ideals.
Even if biquadratic extensions may not be optimal for this algorithm \cite{Klein, Murphy}, a more general form of our results on ideal combination could lead to better computational performance. 

Finally, from a theoretical point of view, it may be worthy to investigate which algebraic properties are preserved, as happens for divisibility, by first-degree prime ideals combination.


\section*{acknowledgments}
The authors want to thank professors Massimiliano Sala and Michele Elia for their useful comments and suggestions.







\end{document}